\documentclass[a4paper]{amsproc}
\usepackage{graphicx,amssymb,amsmath}
\begin{document}

\title{Some tastings in Morales-Ramis Theory}

\author{Primitivo Bel\'en Acosta-Hum\'anez}

\address[P. Acosta-Hum\'anez]{Professor PAC, Instituto Superior de Formaci\'on Docente Salom\'e Ure\~na, Recinto Emilio Prud'Homme, Santiago de los Caballeros -- Dominican Republic. Honorary Researcher, Facultad de Ciencias B\'asicas y Biom\'edicas, Universidad Sim\'on Bol\'{\i}var, Barranquilla -- Colombia}

\email{primitivo.acosta-humanez@isfodosu.edu.do}

\author{Germ\'an Jim\'enez Blanco}

\address[G. Jim\'enez Blanco]{Professor, Department of Mathematics and Statistics,Universidad del Norte, Barranquilla - Colombia}

\email{gjimenez@uninorte.edu.co}
\maketitle
\begin{abstract}
In this paper we present a short material concerning to some results in Morales-Ramis theory, which relates two different notions of integrability: Integrability of Hamiltonian Systems through Liouville Arnold Theorem and Integrability of Linear Differential Equations through Differential Galois Theory. As contribution, we obtain the abelian differential Galois group of the variational equation related to a bi-parametric Hamiltonian system,
\end{abstract}

\section{Introduction}
To understand the Morales-Ramis theory, we need to introduce two different notions of integrability, the
integrability of Hamiltonian systems in Liouville sense and the
integrability of linear differential equations in Picard-Vessiot sense.

\subsection{Integrability of Hamiltonian Systems}

From the previous section, let us consider a $n$ degrees of freedom
hamiltonian $H$. 
 The equations of the flow of the hamiltonian system, in a
system of canonical coordinates, $x_1,\ldots,x_n,y_1,\ldots,y_n$, are
written 
\begin{equation*}
\dot x = \frac{\partial H}{\partial y}, \quad \dot
y = - \frac{\partial H}{\partial x} ,
\end{equation*}
and they are known as \emph{Hamilton equations}. 
We recall that the \emph{Poisson bracket} between $f(x_1,x_2,y_1,y_2)$ and $g(x_1,x_2,y_1,y_2)$ is given by
$$\{f,g\}=\sum_{k=1}^n\left({\partial f\over \partial y_k}{\partial g\over \partial x_k}-{\partial f\over \partial y_k}{\partial g\over \partial x_k}\right).$$ We say that $f$ and $g$ are in involution when $\{f,g\}=0$, also we say in this case that $f$ and $g$ commute under the Poisson bracket.
In this way, we can write the Hamilton equations as follows:
\begin{equation*}
\dot x = \{H, x\}, \quad \dot
y = \{ H, y\}.
\end{equation*}

A hamiltonian $H$ in $\mathbb{C}^{2n}$ is called \emph{integrable in the sense of Liouville} if
there exist $n$ independent first integrals of the hamiltonian system in involution. We will
say that $H$ in integrable \emph{by terms of rational functions} if we can
find a complete set of integrals within the family of rational functions.
Respectively, we can say that $H$ is integrable \emph{by terms of
meromorphic functions} if we can find a complete set of integrals within the
family of meromorphic functions. The following examples can illustrate a
procedure to obtain first integrals, see also \cite{audin} for further
examples.\medskip

\noindent {\bf Example.} 
Consider the following hamiltonian, which we will proof that is integrable. 
\begin{equation}
H\left( q_{1},q_{2},p_{1},p_{2}\right) =\frac{1}{2}p_{1}^{2}+\frac{1}{2}%
p_{2}^{2}-2q_{1}^{3}-6q_{1}q_{2}^{2},  \label{ex2g}
\end{equation}%
the Hamilton equation for $H$ are. 
\begin{eqnarray*}
\overset{\cdot }{q}_{1} &=&p_{1}\qquad \overset{\cdot }{q}_{2}=p_{2} \\
\overset{\cdot }{p}_{1} &=&6q_{1}^{2}+6q_{2}^{2}\qquad \overset{\cdot }{p}%
_{2}=12q_{1}q_{2}.
\end{eqnarray*}%
Setting $q_{1}=x_{1}-x_{2}\qquad q_{2}=x_{1}+x_{2}$, then 
\begin{eqnarray*}
\widehat{H}\left( x_{1},x_{2},p_{1}\left( y_{1},y_{2}\right) ,p_{2}\left(
y_{1},y_{2}\right) \right) &=&\frac{1}{2}p_{1}^{2}\left( y_{1},y_{2}\right) +%
\frac{1}{2}p_{2}^{2}\left( y_{1},y_{2}\right) -2\left( x_{1}-x_{2}\right)
^{3}-6\left( x_{1}-x_{2}\right) \left( x_{1}+x_{2}\right) ^{2} \\
&=&\frac{1}{2}p_{1}^{2}\left( y_{1},y_{2}\right) +\frac{1}{2}p_{2}^{2}\left(
y_{1},y_{2}\right) -8x_{1}^{3}+8x_{2}^{3}.
\end{eqnarray*}%
From the Hamilton equations corresponding to $\widehat{H}$ we have 
\begin{eqnarray*}
\overset{\cdot }{x}_{1} &=&\frac{\partial \widehat{H}}{\partial y_{1}}=p_{1}%
\frac{\partial p_{1}}{\partial y_{1}}+p_{2}\frac{\partial p_{2}}{\partial
y_{1}} \\
\overset{\cdot }{x}_{2} &=&\frac{\partial \widehat{H}}{\partial y_{2}}=p_{1}%
\frac{\partial p_{1}}{\partial y_{2}}+p_{2}\frac{\partial p_{2}}{\partial
y_{2}},
\end{eqnarray*}%
since 
\begin{eqnarray*}
\overset{\cdot }{x}_{1}+\overset{\cdot }{x}_{2} &=&\overset{\cdot }{q}%
_{2}=p_{2} \\
\overset{\cdot }{x}_{1}-\overset{\cdot }{x}_{2} &=&\overset{\cdot }{q}%
_{1}=p_{1}
\end{eqnarray*}%
\begin{eqnarray*}
p_{1}\left( \frac{\partial p_{1}}{\partial y_{1}}+\frac{\partial p_{1}}{%
\partial y_{2}}\right) +p_{2}\left( \frac{\partial p_{2}}{\partial y_{1}}+%
\frac{\partial p_{2}}{\partial y_{2}}\right) &=&p_{2} \\
p_{1}\left( \frac{\partial p_{1}}{\partial y_{1}}-\frac{\partial p_{1}}{%
\partial y_{2}}\right) +p_{2}\left( \frac{\partial p_{2}}{\partial y_{1}}-%
\frac{\partial p_{2}}{\partial y_{2}}\right) &=&p_{1}.
\end{eqnarray*}%
A solution for the above system is 
\begin{eqnarray*}
p_{1} &=&\frac{y_{1}-y_{2}}{2} \\
p_{2} &=&\frac{y_{1}+y_{2}}{2},
\end{eqnarray*}%
then 
\begin{eqnarray*}
\widehat{H}\left( x_{1},x_{2},y_{1},y_{2}\right) &=&\frac{1}{2}\left( \frac{%
y_{1}-y_{2}}{2}\right) ^{2}+\frac{1}{2}\left( \frac{y_{1}+y_{2}}{2}\right)
^{2}-8x_{1}^{3}+8x_{2}^{3} \\
&=&\frac{1}{4}y_{1}^{2}+\frac{1}{4}y_{2}^{2}-8x_{1}^{3}+8x_{2}^{3}.
\end{eqnarray*}%
For instance, a first integral corresponding to $\widehat{H}$ is 
\begin{equation*}
\widehat{I}=\frac{1}{4}y_{1}^{2}-8x_{1}^{3}
\end{equation*}%
and therefore another first integral corresponding to $H$ is 
\begin{equation*}
I=\frac{1}{4}\left( p_{1}+p_{2}\right) ^{2}-\left( q_{1}+q_{2}\right) ^{3}.
\end{equation*}%
In this way, we have that $H$ and $I$ are first integrals, independent and in involution due to 
\begin{eqnarray*}
\frac{d I}{d t} &=&\frac{1}{2}\left( p_{1}+p_{2}\right) \left( 
\overset{\cdot }{p}_{1}+\overset{\cdot }{p}_{2}\right) -3\left(
q_{1}+q_{2}\right) ^{2}\left( \overset{\cdot }{q}_{1}+\overset{\cdot }{q}%
_{2}\right) \\
&=&\frac{1}{2}\left( p_{1}+p_{2}\right) \left(
6q_{1}^{2}+6q_{2}^{2}+12q_{1}q_{2}\right) -3\left( q_{1}+q_{2}\right)
^{2}\left( p_{1}+p_{2}\right) \\
&=&3\left( p_{1}+p_{2}\right) \left( q_{1}+q_{2}\right) ^{2}-3\left(
q_{1}+q_{2}\right) ^{2}\left( p_{1}+p_{2}\right) =0
\end{eqnarray*}
and $\{H, I\}=0$.

We denote by $X_H$ the hamiltonian vector field, that is, the right-hand side of the hamilton equations. In a general way, we deal with non-linear hamiltonian systems. For suitability, without lost of generality, we can consider hamiltonian systems with two degrees of freedom, that is a hamiltonian $H$ in 
$\mathbb{C}^{4}$. Let $\Gamma$ be an integral curve of $X_H$, being $\Gamma$ parametrized by $\gamma\colon t\mapsto (q_1(t),q_2(t),p_1(t),p_2(t))$, the \emph{first variational equation} (VE) along $\Gamma$ is given by, 
\begin{equation*}
\left(%
\begin{array}{c}
\dot \xi_1 \\ 
\dot\xi_2\\
\dot \xi_3 \\ 
\dot\xi_4%
\end{array}%
\right) = \left(%
\begin{array}{cccc}
\frac{\partial^2 H}{\partial p_1\partial q_1}(\gamma(t)) & \frac{\partial^2 H%
}{\partial p_1 \partial q_2}(\gamma(t)) & \frac{\partial^2 H}{\partial p_1^2}(\gamma(t)) & \frac{\partial^2 H%
}{\partial p_1 \partial p_2}(\gamma(t))\\
\frac{\partial^2 H}{\partial p_2\partial q_1}(\gamma(t)) & \frac{\partial^2 H%
}{\partial p_2 \partial q_2}(\gamma(t)) & \frac{\partial^2 H}{\partial p_2\partial p_1}(\gamma(t)) & \frac{\partial^2 H%
}{\partial p_2^2}(\gamma(t))\\
- \frac{\partial^2 H}{\partial q_1^2}(\gamma(t)) & -\frac{%
\partial^2 H}{\partial q_1 \partial q_2}(\gamma(t))&- \frac{\partial^2 H}{\partial q_1 \partial p_1}(\gamma(t)) & -\frac{%
\partial^2 H}{\partial q_1 \partial p_2}(\gamma(t))\\- \frac{\partial^2 H}{\partial q_2 \partial q_1}(\gamma(t)) & -\frac{%
\partial^2 H}{\partial q_2^2}(\gamma(t))&- \frac{\partial^2 H}{\partial q_2 \partial p_1}(\gamma(t)) & -\frac{%
\partial^2 H}{\partial q_2 \partial p_2}(\gamma(t))
\end{array}%
\right) \left(%
\begin{array}{c}
\xi_1 \\ 
\xi_2\\
 \xi_3 \\ 
\xi_4%
\end{array}
\right).
\end{equation*}

\subsection{Picard-Vessiot Theory.}

The Picard-Vessiot theory is the Galois theory of linear differential
equations. In the classical Galois theory, the main object is a group of
permutations of the roots, while in the Picard-Vessiot theory is a linear
algebraic group. In the remainder of this paper we only work, as particular
case, with linear differential equations of second order 
\begin{equation}  \label{LDE}
y^{\prime \prime }+ay^{\prime }+by=0,\quad a,b\in \mathbb{C}(x).  \tag{LDE}
\end{equation}
Suppose that $y_1, y_2$ is a fundamental system of solutions of the
differential equation. This means that $y_1, y_2$ are linearly independent
over $\mathbb{C}$ and every solution is a linear combination of these two.
Let $L = \mathbb{C}(x)\langle y_1, y_2 \rangle = \mathbb{C}(x)(y_1, y_2,
y_1^{\prime }, y_2^{\prime })$, that is the smallest differential field
containing to $\mathbb{C}(x)$ and $\{y_{1},y_{2}\}.$

The group of all differential automorphisms of $L$ over $\mathbb{C}(x)$ is
called the {\emph{Galois group}} of $L$ over $\mathbb{C}(x)$ and denoted by $%
Gal(L/\mathbb{C}(x))$ or also by $Gal^L_{\mathbb{C}(x)}$. This means that
for $\sigma\colon L\to L$, $\sigma(a^{\prime })=\sigma^{\prime }(a)$ and $%
\forall a\in \mathbb{C}(x),$ $\sigma(a)=a$.

If $\sigma \in Gal(L/\mathbb{C}(x))$ then $\sigma y_1, \sigma y_2$ is
another fundamental system of solutions of the linear differential equation.
Hence there exists a matrix $A= 
\begin{pmatrix}
a & b \\ 
c & d%
\end{pmatrix}
\in GL(2,\mathbb{C})$, such that 
\begin{equation*}
\sigma 
\begin{pmatrix}
y_{1} \\ 
y_{2}%
\end{pmatrix}
= 
\begin{pmatrix}
\sigma y_{1} \\ 
\sigma y_{2}%
\end{pmatrix}
=A 
\begin{pmatrix}
y_{1} \\ 
y_{2}%
\end{pmatrix}
.
\end{equation*}

{\bf Theorem 1.}
The Galois group $G=Gal(L/\mathbb{C}(x))$ is an algebraic subgroup of $GL(2,%
\mathbb{C})$. Moreover, the Galois group of a reduced linear differential
equation 
\begin{equation}  \label{LDE}
\xi^{\prime \prime }=r\xi,\quad r\in \mathbb{C}(x),  \tag{RLDE}
\end{equation}
is an algebraic subgroup $SL(2,\mathbb{C})$.

Let $F\subset L$ be a differential field extension, and let $\eta$ be and
element of $L$, then

\begin{enumerate}
\item $\eta$ is {\emph{algebraic}} over $F$ if $\eta$ satisfies a polynomial
equation with coefficients in $F$, i.e. $\eta$ is an algebraic function of
one variable.

\item $\eta$ is {\emph{primitive}} over $F$ if $\eta^{\prime }\in F$, i.e. $%
\eta = \int f$ for some $f \in F$.

\item $\eta$ is {\emph{exponential}} over $F$ if $\eta^{\prime }/\eta \in F$%
, i.e. $\eta = e^{\int f}$ for some $f \in F$.
\end{enumerate}

Algebraic, primitive and exponential functions are called Liouvillian
functions. Thus a Liouvillian function is built up using algebraic
functions, integrals and exponentials. In the case $F=\mathbb{C}(x)$ we get,
for instance logarithmic, trigonometric functions, but not special functions
such that the Airy functions. Solvability of LDE is expressed in terms of
these function. This is one of the main results of Picard-Vessiot theory.

{\bf Theorem 2.}
A linear differential equation is solvable integrable by terms of,
Liouvillian functions, if and only if the connected component of the
identity element of its Galois group is a solvable group.

The Cauchy-Euler equation with Galois group in $SL(2,\mathbb{C})$ is given by
$$y''={m(m+1)\over x^2}.$$ The integrability and Galois group of this equation has been studied by the first author in his PhD thesis, where is proved that the differential Galois group is abelian for all $m\in \mathbb{C}$.

\section{Morales Ramis Theory}

We want to relate integrability of hamiltonian systems to Picard-Vessiot
theory. The following theorems treat this problem.

{\bf Theorem 3 (Morales-Ramis, \cite{mora1}).}
Let $H$ be a Hamiltonian in $\mathbb{C}^{2n}$, and $\gamma$ a
particular solution such that the \textrm{NVE} has regular (resp. irregular)
singularities at the points of $\gamma$ at infinity. Then, if $H$ is
completely integrable by terms of meromorphic (resp. rational) functions,
then the Identity component of Galois Group of the \textrm{NVE} is abelian.

To understand completely this technical result, it is required a formal study of concerning to differential Galois theory and Morales-Ramis theory. We can illustrat this theorem through the following examples, but only for a basic level.

{\bf Example.}
We study the conditions of the previous theorem in the hamiltonian system given in equation \eqref{ex2g}, which is given by

$$H\left( q_{1},q_{2},p_{1},p_{2}\right) =\frac{1}{2}p_{1}^{2}+\frac{1}{2}%
p_{2}^{2}-2q_{1}^{3}-6q_{1}q_{2}^{2}.$$
The hamilton equations are:%
\begin{eqnarray*}
\overset{\cdot }{q_{1}} &=&p_{1},\qquad \overset{\cdot }{q_{2}}=p_{2} \\
\overset{\cdot }{p_{1}} &=&6q_{1}^{2}+6q_{1}^{2} \\
\overset{\cdot }{p_{2}} &=&12q_{1}q_{2}
\end{eqnarray*}%
taking the invariant plane $q_{2}=p_{2}=0.$ we have $\overset{\cdot \cdot }{%
q_{1}}=6q_{1}$ a solution for this equation is $q_{1}\left( t\right) =\frac{1%
}{t^{2}}$, and the variational equation is:%
\begin{equation*}
\left[ 
\begin{array}{cccc}
0 & 0 & 1 & 0 \\ 
0 & 0 & 0 & 1 \\ 
\frac{12}{t^{2}} & 0 & 0 & 0 \\ 
0 & \frac{12}{t^{2}} & 0 & 0%
\end{array}%
\right] \left[ 
\begin{array}{c}
\xi _{1} \\ 
\xi _{2} \\ 
\xi _{3} \\ 
\xi _{4}%
\end{array}%
\right] =\left[ 
\begin{array}{c}
\overset{\cdot }{\xi _{1}} \\ 
\overset{\cdot }{\xi _{2}} \\ 
\overset{\cdot }{\xi _{3}} \\ 
\overset{\cdot }{\xi _{4}}%
\end{array}%
\right]
\end{equation*}%
then $12\xi _{1}=t^{2}\overset{\cdot \cdot }{\xi _{1}}$ , which corresponds to a 
Cauchy-Euler equation, thus, the Galois group is abelian due to the hamiltonian system is
integrable.

The following examples were taken from \cite{acbl}.\medskip

{\bf Example.}
Consider the hamiltonian 
\begin{equation}
H=\frac{p_{1}^{2}+p_{2}^{2}}{2}-Q(q_{1})\frac{q_{2}^{2}}{2}+\beta
(q_{1},q_{2})q_{2}^{3},  \label{ex3g}
\end{equation}
where $Q(q_{1})$ is a polynomial and $\beta \left( q_{1},q_{2}\right) $ is a
function of two variables with continuous partial derivative and $%
\lim\limits_{q_{2}\rightarrow 0}\frac{\partial ^{j}\beta \left(
q_{1},q_{2}\right) }{\partial q_{2}^{j}}<\infty \ \ 0\leq j\leq 2.$\newline
The hamilton equations are:%
\begin{eqnarray*}
\overset{\cdot }{q}_{1} &=&p_{1},\qquad \overset{\cdot }{q}_{2}=p_{2} \\
\overset{\cdot }{p}_{1} &=&Q^{^{\prime }}(q_{1})\frac{q_{2}^{2}}{2}-\frac{%
\partial \beta (q_{1},q_{2})}{\partial q_{1}}q_{2}^{3} \\
\overset{\cdot }{p}_{2} &=&Q\left( q_{1}\right) q_{2}-\frac{\partial \beta
(q_{1},q_{2})}{\partial q_{2}}q_{2}^{3}-3\beta \left( q_{1},q_{2}\right)
\left( q_{2}\right) ^{2}
\end{eqnarray*}%
taking the invariant plane $q_{2}=p_{2}=0.$ we have $q_{1}\left( t\right)
=at+b$ and the variational equation is:%
\begin{equation*}
\left[ 
\begin{array}{cccc}
0 & 0 & 1 & 0 \\ 
0 & 0 & 0 & 1 \\ 
0 & 0 & 0 & 0 \\ 
0 & Q\left( q_{1}\right)  & 0 & 0%
\end{array}%
\right] \left[ 
\begin{array}{c}
\xi _{1} \\ 
\xi _{2} \\ 
\xi _{3} \\ 
\xi _{4}%
\end{array}%
\right] =\left[ 
\begin{array}{c}
\overset{\cdot }{\xi _{1}} \\ 
\overset{\cdot }{\xi _{2}} \\ 
\overset{\cdot }{\xi _{3}} \\ 
\overset{\cdot }{\xi _{4}}%
\end{array}%
\right] 
\end{equation*}%
then $\xi _{3}=\overset{\cdot }{\xi _{1}},\qquad \xi _{4}=\overset{\cdot }{%
\xi _{2}},\qquad Q\left( q_{1}\right) \xi _{2}=\overset{\cdot }{\xi _{4}}$
hence $Q\left( q_{1}\right) \xi _{2}=\overset{\cdot \cdot }{\xi_2}.$ \newline
If $Q\left( q_{1}\right) $ is a polynomial then the Galois group is
not abelian, although in some case is solvable (see \cite{acbl}), hence the Hamiltonian System
is not integrable by Morales-Ramis Theorem. 

{\bf Example.}
Consider the following hamiltonian 
\begin{equation}
H=\frac{p_{1}^{2}+p_{2}^{2}}{2}-\frac{\lambda _{4}}{(\lambda _{2}+2\lambda
_{3}q_{1})^{2}}+\lambda _{0}-\lambda _{1}q_{2}^{2}-\lambda
_{2}q_{1}q_{2}^{2}-\lambda _{3}q_{1}^{2}q_{2}^{2}+\beta
(q_{1},q_{2})q_{2}^{3},
\end{equation}%
Where $\lambda _{i}\in 
\mathbb{C}
,\ \ $with $\lambda _{3}\neq 0$, \  $\beta \left( q_{1},q_{2}\right) $ is a
function of two variables with continuous partial derivative and $%
\lim\limits_{q_{2}\rightarrow 0}\frac{\partial ^{j}\beta \left(
q_{1},q_{2}\right) }{\partial q_{2}^{j}}<\infty \ \ 0\leq j\leq 2.$

The hamilton equations are:%
\begin{eqnarray*}
\overset{\cdot }{q}_{1} &=&p_{1},\qquad \overset{\cdot }{q}_{2}=p_{2} \\
\overset{\cdot }{p}_{1} &=&\frac{-4\lambda _{3}\lambda _{4}}{(\lambda
_{2}+2\lambda _{3}q_{1})^{3}}+\lambda _{2}q_{2}^{2}+2\lambda
_{3}q_{1}q_{2}^{2}-\frac{\partial \beta }{\partial q_{1}}\left(
q_{1},q_{2}\right) q_{2}^{3} \\
\overset{\cdot }{p_{2}} &=&2\lambda _{1}q_{2}+2\lambda
_{2}q_{1}q_{2}+2\lambda _{3}q_{1}^{2}q_{2}-\frac{\partial \beta }{\partial
q_{2}}\left( q_{1},q_{2}\right) q_{2}^{3}-3\beta (q_{1},q_{2})q_{2}^{2}.
\end{eqnarray*}%
Taking $\ q_{2}=p_{2}=0$ and setting $H\left( q_{1},0,p_{1},0\right) =h$, we see that
\begin{eqnarray*}
h &=&\frac{1}{2}p_{1}^{2}-\dfrac{\lambda _{4}}{(\lambda _{2}+2\lambda
_{3}q_{1})^{2}}+\lambda _{0} \\
\overset{\cdot }{q_{1}} &=&p_{1}=\sqrt{2h+\dfrac{2\lambda _{4}}{(\lambda
_{2}+2\lambda _{3}q_{1})^{2}}-2\lambda _{0}}.
\end{eqnarray*}%
Now, we pick $h=\lambda _{0}$, thus we have 
\begin{eqnarray*}
\frac{dq_{1}}{dt} &=&\sqrt{\dfrac{\lambda _{4}}{(\lambda _{2}+2\lambda
_{3}q_{1})^{2}}} \\
\lambda _{2}q_{1}+\lambda _{3}q_{1}^{2} &=&\pm \sqrt{\lambda _{4}}t+c.
\end{eqnarray*}%
For instance the variational equation is%
\begin{equation*}
\left[ 
\begin{array}{cccc}
0 & 0 & 1 & 0 \\ 
0 & 0 & 0 & 1 \\ 
\frac{24\lambda _{3}^{2}\lambda _{4}}{(\lambda _{2}+2\lambda _{3}q_{1})^{4}}
& 0 & 0 & 0 \\ 
0 & 2\lambda _{1}+2\lambda _{2}q_{1}+2\lambda _{3}q_{1}^{2} & 0 & 0%
\end{array}%
\right] \left[ 
\begin{array}{c}
\xi _{1} \\ 
\xi _{2} \\ 
\xi _{3} \\ 
\xi _{4}%
\end{array}%
\right] =\left[ 
\begin{array}{c}
\overset{\cdot }{\xi _{1}} \\ 
\overset{\cdot }{\xi _{2}} \\ 
\overset{\cdot }{\xi _{3}} \\ 
\overset{\cdot }{\xi _{4}}%
\end{array}%
\right] 
\end{equation*}%
then 
\begin{equation}
\left( 2\lambda _{1}+2\lambda _{2}q_{1}+2\lambda _{3}q_{1}^{2}\right) \xi
_{2}=\overset{\cdot \cdot }{\xi _{2}}  \label{L13}
\end{equation}%
replacing in equation $\left( \ref{L13}\right) $ we have 
\begin{equation}
p\left( t\right) \xi _{2}=\overset{\cdot \cdot }{\xi _{2}}  \label{L14}
\end{equation}%
where $p\left( t\right) =2\lambda _{1}+2\left( \pm \sqrt{\lambda _{4}}%
t+c\right) .$

and consequently the Galois group of $\left( \ref{L14}\right) $ is
not abelian hence the hamiltonian system is not integrable. 

{\bf Example.}
Consider the following hamiltonian%
\begin{equation}
H=\frac{1}{2}p_{1}^{2}+\frac{1}{2}p_{1}^{2}+\frac{1}{%
aq_{1}^{3}+bq_{1}^{2}q_{2}+cq_{2}^{3}}
\end{equation}%
Where $a,b,c\in 
\mathbb{C}
.$ The hamilton equations are:%
\begin{eqnarray*}
\overset{\cdot }{q}_{1} &=&p_{1} \\
\overset{\cdot }{q}_{2} &=&p_{2} \\
\overset{\cdot }{p}_{1} &=&\frac{3aq_{1}^{2}+2bq_{1}q_{2}}{\left(
aq_{1}^{3}+bq_{1}^{2}q_{2}+cq_{2}^{3}\right) ^{2}} \\
\overset{\cdot }{p_{2}} &=&\frac{bq_{1}^{2}+3cq_{2}^{2}}{\left(
aq_{1}^{3}+bq_{1}^{2}q_{2}+cq_{2}^{3}\right) ^{2}}
\end{eqnarray*}%
taking the invariant plane $q_{1}=p_{1}=0$ we have 
\begin{eqnarray*}
\overset{\cdot }{p_{2}} &=&\frac{3}{cq_{2}^{4}}\qquad \overset{\cdot \cdot }{%
q_{2}}=\frac{3}{cq_{2}^{4}} \\
q_{2}\left( t\right) &=&\left( \dfrac{-25}{2c}\right) ^{\dfrac{1}{5}}t^{%
\dfrac{2}{5}}
\end{eqnarray*}%
$Z\left( t\right) =\left( 0,\left( \dfrac{-25}{2c}\right) ^{\dfrac{1}{5}}t^{%
\dfrac{2}{5}},0,\dfrac{-5}{c\left( -\frac{25}{2c}\right) ^{\frac{4}{5}}t^{%
\frac{3}{5}}}\right) .\ \ $the variational equation is 
\begin{equation*}
\left[ 
\begin{array}{c}
\overset{\cdot }{\xi _{1}} \\ 
\overset{\cdot }{\xi _{2}} \\ 
\overset{\cdot }{\xi _{3}} \\ 
\overset{\cdot }{\xi _{4}}%
\end{array}%
\right] =\left[ 
\begin{array}{cccc}
0 & 0 & 1 & 0 \\ 
0 & 0 & 0 & 1 \\ 
\frac{-4b}{25t^{2}} & 0 & 0 & 0 \\ 
0 & \frac{24}{-25t^{2}} & 0 & 0%
\end{array}%
\right] \left[ 
\begin{array}{c}
\xi _{1} \\ 
\xi _{2} \\ 
\xi _{3} \\ 
\xi _{4}%
\end{array}%
\right]
\end{equation*}%
$\overset{\cdot }{q}_{1}=p_{1}\qquad \overset{\cdot }{p}_{1}=\frac{-4b}{%
25t^{2}}q_{1}$ then $\overset{\cdot \cdot }{q}_{1}=\frac{-4b}{25t^{2}}q_{1}$%
is a Cauchy-Euler equation, the Galois group is the identity, this group is
abelian but we cannot state that the Hamiltonian system is integrable.

The previous theorem of Morales-Ramis was extended by Morales-Ramis-Sim\'o.
 
{\bf Theorem 4. (Morales-Ramis-Simo, \cite{morasi}).}
Let $H$ be a Hamiltonian in $\mathbb{C}^{2n}$, and $\gamma $ a
particular solution such that the \textrm{NVE} has regular (resp. irregular)
singularities at the points of $\gamma $ at infinity. Then, if $H$ is
completely integrable by terms of meromorphic (resp. rational) functions,
then the identity component of Galois Group of any linearised high order
variational equation is abelian.

The following example illustrate the way to compute the second order variational equation.

Consider the Hamiltonian

\begin{equation}
H=\frac{1}{2}p_1^2+\frac{1}{2}p_2^2+\frac{1}{2}a_1q_1^2+\frac{1}{2}a_2q_2^2+%
\frac{1}{4}a_5q_1^4+\frac{1}{4}a_3q_2^4+\frac{1}{2}a_4q_1^2q_2^2,
\end{equation}
where the Hamilton equations are given by 
\begin{equation}
\begin{array}{lll}
\dot{q}_1 & = & p_1 \\ 
\dot{q}_2 & = & p_2 \\ 
\dot{p}_1 & = & -a_1q_1-a_5q_1^3-a_4q_1q_2^2 \\ 
\dot{p}_2 & = & -a_2q_2-a_3q_2^3-a_4q_1^2q_2.%
\end{array}%
\end{equation}
Taking as invariant plane $\Gamma=\{(q_1,q_2,p_1,p_2): q_2=p_2=0\}$ we
obtain first variational equation

\begin{equation}
\dot \xi^{(1)}=A(t)\xi^{(1)}, \quad A(t)= 
\begin{pmatrix}
0 & 0 & 1 & 0 \\ 
0 & 0 & 0 & 1 \\ 
-a_1-3a_5q_1^2 & 0 & 0 & 0 \\ 
0 & -a_2-a_4q_1^2 & 0 & 0%
\end{pmatrix}%
,
\end{equation}
where $\xi^{(1)}=(\xi^{(1)}_1,\xi^{(1)}_2,\xi^{(1)}_3,\xi^{(1)}_4)^T$ and $%
q_1=q_1(t)$, being $(q_1(t),0,\dot q_1(t),0)$ a particular solution of the
hamiltonian system over the invariant plane.

The second variational equation is given by

\begin{equation}
\dot \xi^{(2)}=A(t)\xi^{(2)}+f(t), \quad f(t)=%
\begin{pmatrix}
0 \\ 
0 \\ 
-3a_5(\xi^{(1)}_1)^2-a_4q_1(\xi^{(1)}_2)^2 \\ 
-2a_4q_1\xi^{(1)}_1\xi^{(1)}_2,%
\end{pmatrix}%
\end{equation}
where $\xi^{(2)}=(\xi^{(2)}_1,\xi^{(2)}_2,\xi^{(2)}_3,\xi^{(3)}_4)^T$.

The following is our original contribution to this paper.

{\bf{Proposition.}} 
Assume the Hamiltonian system given by $$H={p_1^2+p_2^2\over 2}-{1 \over aq_1^m+bq_2^m}, \quad a\neq 0, \quad m>2.$$ The differential Galois group of the variational equation corresponding to the invariant plane $q_2=p_2=0$ and energy level $h=0$, is virtually abelian. Furthermore, the Galois group is independent of the choice of $a$ and $b$.  

\noindent\textit{Proof.}
The subsystem in the invariant plane is $$h={p_1^2\over 2}-{1 \over aq_1^m}, \quad a\neq 0, \quad m>2,$$ then we obtain a particular solution for $q_1$ given by $$q_1(t)=\left({m+2\over \sqrt{2a}t}\right)^{2\over m+2},$$ for instance the variational equation is given by \begin{equation}
\dot \xi=A(t)\xi, \quad A(t)= 
\begin{pmatrix}
0 & 0 & 1 & 0 \\ 
0 & 0 & 0 & 1 \\ 
{m(m+1)\over aq_1^{m+2}} & 0 & 0 & 0 \\ 
0 & 0 & 0 & 0
\end{pmatrix}, \quad \xi=(\xi_1,\xi_2,\xi_3,\xi_4)^T.
\end{equation}
Thus, we arrive to the Cauchy-Euler equation $${\frac {d^{2}}{d{t}^{2}}}\xi_1 ={\frac {2m
 \left( m+1 \right) }{ \left( m+2
 \right) ^{2}{t}^{2}}}\xi_1
$$
and for instance the Galois group is always abelian.


\end{document}